\documentclass[12pt]{article}
\usepackage{amsmath}
\usepackage{amssymb}
\usepackage{graphicx}
\date{}

\begin{document}

\title{EMD Algorithm with Approximate Zero Crossings}

\author{ Mayer Humi\\
Department of Mathematical Sciences, \\
Worcester Polytechnic Institute, \\
Worcester, MA 01609, USA  \thanks {e-mail: mhumi@wpi.edu. } 
}

\maketitle

\begin{abstract}

The classical EMD algorithm has been used extensively in the literature to 
decompose signals that contain nonlinear waves. However when a signal contain 
two or more frequencies that are close to one another the decomposition might 
fail. In this paper we propose a new formulation of this algorithm which is 
based on the zero crossings of the signal and show that it performs well 
even when the classical algorithm fail. We address also
the filtering properties and convergence rate of the new algorithm versus 
the classical EMD algorithm.  

\end{abstract}

\thispagestyle{empty}

\vspace{2em}

\newpage

\section{Introduction}

In scientific literature there exist many classical sets of functions 
which can decompose a signal in terms of "simple" functions. For 
example Taylor or Fourier expansions are used routinely in scientific and 
engineering applications (and many other exist). However in all these
expansions the underlying functions are not intrinsic to the signal itself
and a precise approximation to the original signal might require a large
number of terms. This problem become even more acute when the signal is 
non-stationary and the process it represents is nonlinear.

To overcome this problem many researchers used in the past 
the "principal component algorithm" (PCA) to come up with an
"adaptive" set of functions which approximate a given signal. 
A new approach to this problem emerged in the late 1990's when a NASA team 
has developed the "Empirical Mode Decomposition" algorithm(EMD) which 
attempts to decompose a signal in terms of it "intrinsic mode functions"
(IMF) through "sifting algorithm". A patent for this algorithm has been 
issued [1]. 

The EMD algorithm is based on the following quote [2]:
"According to Drazin the first step of data analysis is to examine
the data by eye. From this examination, one can immediately identify
the different scales directly in two ways: by the time lapse between
successive alterations of local maxima and minima and by the time lapse
between the successive zero crossings....We have decided to to adopt the
time lapse between successive extrema as the definition of the time
scale for the intrinsic oscillatory mode" 

A step by step description of the EMD sifting algorithm is as follows:
\begin{enumerate}
\item Let be given a function $f(t)$ which is sampled at discrete times 
$\{t_k, k=1,\ldots n\}$.
\item let $h_0(k)=f(t_k)$.
\item Identify the max and min of $h_0(k)$.
\item Create the cubic spline curve $M_x$ that connects the maxima points.
Do the same for the minima $M_n$. This creates an envelope for
$h_0(k)$.
\item At each time $t_k$ evaluate the mean $m_k$ of $M_x$ and $M_n$
($m_k$ is referred to as the sifting function).
\item Evaluate $h_1(k)=h_0(k)-m_k$.
\item If the norm of $||h_0-h_1|| < \epsilon$ for some predetermined 
$\epsilon$ set the first intrinsic mode function $IMF_1=h_1$ (and stop).
\item If the criteria of (7) are not satisfied  set $h_0(k)=h_1(k)$
and return to (3) ("Sifting process"). 
\end{enumerate}

The algorithm has been applied successfully in various physical 
applications [1-6]. However as has been observed by Flandrin [3] and others
the EMD algorithm fails in many cases where the data contains two
or more frequencies which are close to each other.

To overcome this difficulty we propose hereby a modification of the EMD
algorithm by replacing steps $4$ and $5$ in the description above by 
the following:

4. find the midpoints between two consecutive maxima and minima
and let $N_k$ be the values of $h_0$ at these points. \\
5. Create the spline curve $m_k$ that connects the points $N_k$.

The essence of this modification is the replacement of the mean 
which is evaluated by the EMD algorithm as the average of the max-min
envelopes by the spline curve of the mid-points between the maxima and minima.
This is in line with the observation by Drazin (which was referred to above)
that the scales inherent to the data can be educed either from the max-min 
or its zero crossing. In the algorithm we propose hereby we mimic the 
"zero-crossings" by the mid-points between the max-min. 
 
It is our objective in this paper to justify this modification of the EMD 
algorithm through some theoretical work and case studies. The plan of the
paper is as follows: In Sec $2$ we provide theoretical justification for the
new algorithm proving that it acts as a high pass filter for certain 
classes of signals. In Sec. $3$ we provide examples of a signal 
composed of two or three close frequencies (with and without noise) where 
the classical EMD algorithm fails but the modified one yields satisfactory 
result. In Sec. $4$ we discuss the convergence rate, resolution and related 
issues concerning the classical and new "midpoint algorithm". 
We end up with some conclusions in Sec $5$.

\setcounter{equation}{0}
\section{Theoretical Justification}

In this section we provide a theoretical justification for the proposed
modified EMD algorithm by analyzing it performance on several generic 
signals which contain several close frequencies. However in this analysis
linear,quadratic and cubic interpolating polynomials will be applied to
represent the midpoints interpolating function (instead of splines). 
To motivate this "replacement" we observe that the coefficients of each 
spline polynomial depend non-locally on the data i.e. these coefficients 
might change if additional data is added. On the other hand Lagrange 
interpolating polynomials depend only the local data.

{\bf Lemma\,1}: Consider a signal of the form
\begin{equation}
\label{0.1}
f(t)= f_1(t)+f_2(t)+f_3(t)
\end{equation}
where
\begin{equation}
\label{0.2}
f_1(t)=\cos(\omega t),\,\,\,f_2(t)=\cos((1+a\epsilon)\omega t),\,\,\,\,
f_3(t)=\cos((1+b\epsilon)\omega t),\,\,\,\, b > a> 0,\,\,\,0 < \epsilon \ll 1.
\end{equation}
Let the projection of the midpoint {\bf linear interpolating function} on 
$f_1(t),\,f_2(t),\, f_3(t)$ over an interval containing five midpoints 
be denoted respectively by $P_{11}$, $P_{21}$, $P_{31}$ then 
$$
P_{i1} = \frac{8\pi}{3}\frac{(b^2+a^2-ba)}{\omega}\epsilon +O(\epsilon^2)
$$
and
\begin{equation}
\label{0.3}
P_{11}-P_{31}= \frac{248\pi}{3}\frac{(b^2+a^2-ba)b}{\omega}\epsilon^3+O(\epsilon^4),
\end{equation}
\begin{equation}
\label{0.4}
P_{11}-P_{21}= \frac{248\pi}{3}\frac{(b^2+a^2-ba)a}{\omega}\epsilon^3+O(\epsilon^4),
\end{equation}
\begin{equation}
\label{0.5}
P_{21}-P_{31}= \frac{248\pi}{3}\frac{(b^2+a^2-ba)(b-a)}{\omega}
\epsilon^3+O(\epsilon^4).
\end{equation}

{\bf Proof}: As a first step we find the approximate location of the 
extrema of $f(t)$ on the interval $[0,6\pi]$. To do so we 
differentiate $f(t)$ and observe that due to the fact that $\epsilon \ll 1$ 
the locations of these points are close to $\frac{n\pi}{\omega}$. Setting 
$t=\frac{n\pi}{\omega}+\eta$ we expand $f'(t)$ in a Taylor series to order $2$
in $\eta$ around $\frac{n\pi}{\omega}$. Then we solve for $\eta$ to obtain 
the approximate locations of the extrema points. Taking the five midpoints
$\{t_1,\ldots,t_5\}$ between these extrema and evaluating 
$f(t_i),\,i=1,\ldots,5$ we construct the linear interpolating 
function $g_1(t)$ between these points. The projection
of $f_i(t)$ on $g_1(t)$ is
\begin{equation}
\label{0.6}
P_{i1}= \displaystyle\int_{t_1}^{t_5} f_i(t)g_1(t) dt.
\end{equation}
Expanding $P_{i1}-P_{j1}$ in a Taylor series in $\epsilon$ one obtains 
(\ref{0.3})-(\ref{0.5}).

{\bf Lemma\,2}: With the same settings as in $Lemma\, 1$ 
let the function $g_2(t)$ consists of the two quadratic polynomials
interpolating $\{t_1,t_2,t_3\}$ and $\{t_3,t_4,t_5\}$ respectively.
The differences between the projections 
\begin{equation}
\label{0.8}
P_{i2}= \displaystyle\int_{t_1}^{t_5} f_i(t)g_2(t) dt
\end{equation}
are given by eqs (\ref{0.3}),(\ref{0.4}) and (\ref{0.5}) respectively
(where $P_{i1}$ is replaced by $P_{i2}$)

{\bf Proof}: As in $Lemma\, 1$  we compute the projections
$P_{i2}$ and expand the results in a Taylor series in $\epsilon$ to 
obtain (\ref{0.3})-(\ref{0.5}) 

{\bf Lemma\,3}: With the same settings as in $Lemma\, 1$ if the projection 
of $f_1(t),\,f_2(t),\, f_3(t)$ is made on the cubic polynomial $g_3(t)$ 
interpolating $\{t_1,t_2,t_3,t_4\}$ then the differences between
the projections
\begin{equation}
\label{0.9}
P_{i3}= \displaystyle\int_{t_1}^{t_4} f_i(t)g_3(t) dt
\end{equation}
are given by 
\begin{equation}
\label{0.10}
P_{13}-P_{33}= \frac{8}{9}\frac{b(b^2+a^2-ba)(49\pi^2-57)}{\omega\pi}\epsilon^3+O(\epsilon^4),
\end{equation}
\begin{equation}
\label{0.11}
P_{13}-P_{23}= \frac{8}{9}\frac{a(b^2+a^2-ba)(49\pi^2-57)}{\omega\pi}\epsilon^3+O(\epsilon^4),
\end{equation}
\begin{equation}
\label{0.12}
P_{23}-P_{33}= \frac{8}{9}\frac{(b-a)(b^2+a^2-ba)(49\pi^2-57)}{\omega\pi}
\epsilon^3+O(\epsilon^4).
\end{equation}

{\bf Theorem\, 1}: As a result of one iteration of the midpoint EMD algorithm
with linear, quadratic or cubic interpolating functions the change in the
projections of the functions $f_i(t) \,i=1,2,3$ on the signal in the interval
$[t_1,t_5])$ ($[t_1,t_4]$ in the cubic case) satisfy 
$\Delta A_{1,j}> \Delta A_{2,j} > \Delta A_{3,j},\,\, j=1,2,3$.
(Here $j$ represents the different interpolating functions).

{\bf Proof}: The projection of $f_i(t)$ on the original signal is
$$
A_{i,j}^0= \displaystyle\int_{t_1}^{t_5} f(t)f_i(t)dt.
$$
After one iteration the signal is represented by 
$$
f^1(t)=f(t)-g_j(t)
$$
and the projection of $f_i(t)$ on $f^1(t)$ is
$$
A_{i,j}^1= \displaystyle\int_{t_1}^{t_5} (f(t)-g_j(t))f_i(t)dt=
A_{i,j}^0-P_{i,j}.
$$
Hence 
$$
\Delta A_{i,j}=A_{i,j}^0-A_{i,j}^1=P_{i,j} > 0
$$
From the results of $lemmas \,1,2,3$ we have that For $j=1,2,3$
$P_{1j}-P_{2j} >0$, $P_{1j}-P_{3j} > 0$ and $P_{2j}-P_{3j} > 0$.
It follows then that $\Delta A_{1,j} > \Delta A_{2,j} > \Delta A_{3,j} > 0$.

We conclude therefore that in the new signal (after one iteration)
the amplitude of $f_3(t)$ will be larger than those of $f_2(t)$ and $f_1(t)$. 
In other words the midpoint EMD algorithm acts as a high pass filter.

We consider now a signal with two close frequencies where a phase shifts 
exists between these two frequencies.  

{\bf Lemma\,4}: Consider a signal of the form
\begin{equation}
\label{0.13}
f(t)= f_4(t)+f_5(t)
\end{equation}
where
\begin{equation}
\label{0.14}
f_4(t)=\cos(\omega t),\,\,\,f_5(t)=\cos[(1+a\epsilon)\omega t+\phi].
\end{equation}
where $a> 0$, and $0 < \epsilon,\phi \ll 1$.
With same setting as in $Lemma \,1$ let the projection of the midpoint 
linear interpolating function (for $f(t)$ defined in (\ref{0.14})) 
on $f_4(t),\,f_5(t)$ over an interval containing five midpoints $\{t_1,\ldots,t_5\}$ be denoted respectively by $P_{41}$, $P_{51}$ then
\begin{equation}
\label{0.15}
P_{41}-P_{51}=\frac{2a\epsilon[6a^2\pi^2\epsilon^2+(5\pi a\epsilon+2\phi)^2]}
{\omega\pi}+O(\epsilon^4,\phi^4)
\end{equation}

{\bf Lemma\,5}:With the same settings as in $Lemma\, 4$ let the function 
$g_5(t)$ consists of the two quadratic polynomials interpolating 
$\{t_1, t_2, t_3\}$ and $\{t_3, t_4, t_5\}$ respectively. The difference 
between the projections $P_{42}$, $P_{52}$ of $f_4(t),\,f_5(t)$ on $g_5(t)$ 
is given by (\ref{0.15}) (where $P_{i1}$ is replaced by $P_{i2}$).

{\bf Lemma\,6}: With the same settings as in $Lemma\,4$ let the function 
$g_6(t)$ consists of the cubic polynomial interpolating 
$\{t_1, t_2, t_3,t_4\}$. The difference between the projections
$P_{43}$, $P_{53}$ of $f_4(t),\,f_5(t)$ on $g_6(t)$ is
\begin{equation}
\label{0.16}
P_{43}-P_{53}=\frac{2a}{3\omega\pi^3}
\left[a^2\pi^2(49\pi^2-57)\epsilon^3+8a\pi\phi(5\pi^2-6)\epsilon^2+2\phi^2(5\pi^2-6)\epsilon\right]+O(\epsilon^4,\phi^4)
\end{equation}

{\bf Theorem\, 2}: As a result of one iteration of the midpoint EMD algorithm
with linear, quadratic or cubic interpolating functions the amplitudes
$B_{i,j},\,i=1,2\,\,j=1,2,3$ (where the index $j$ represents the different
interpolating functions) of the two frequencies present in the
signal (\ref{0.13})-\ref{0.14})) will satisfy
$B_{1,j}< B_{2,j}$. In other words the midpoint EMD algorithm
for this signal is a high pass filter.

{\bf Proof}: The proof is similar to the proof of theorem $1$.

\subsection{Perturbation Analysis}

To investigate the performance of the EMD algorithm (classical and midpoint) 
in the presence of a perturbation (viz. noise) we considered a signal of 
the form
\begin{equation}
\label{5.1}
S_0(t)=\cos(\omega t) +\epsilon f(t)
\end{equation}
where $0 < \epsilon \ll 1$. To analyze this signal we assume that the presence 
of noise (represented by $\epsilon\,f(t)$) does not change (appreciably) 
the location of the extrema in the signal i.e. 
the maximum and minimum are located respectively at the following times
\begin{equation}
\label{5.2}
p_{k}=\frac{2k\pi}{\omega},\,\,\,\, q_k=\frac{(2k+1)\pi}{\omega},\,\,\,\, k=0,1,\ldots.
\end{equation}
The value of the signal at these points is
\begin{equation}
\label{5.3}
S_0(p_k)=1+\epsilon f(\frac{2k\pi}{\omega}),\,\,\,\, 
S_0(q_k)=-1+\epsilon f(\frac{(2k+1)\pi}{\omega}),\,\,\,\, k=0,1,\ldots.
\end{equation}
To apply the classical EMD algorithm to this data one has to compute the spline curves 
$S_{max}(t)$ and $S_{min}(t)$ for the points $(p_k,S(p_k))$ and $(q_k,S(q_k))$ respectively.
The new signal after one iteration of the (classical) EMD algorithm is 
given by
\begin{equation}
\label{5.4}
S_1^c(t)=S_0(t) - \frac{S_{max}(t)+S_{min}(t)}{2} 
\end{equation}
Similarly for the new EMD algorithm we take the midpoints
$d_j=\frac{(2j+1)\pi}{2\omega}$, $j=0,1,\ldots$ between the extrema of the signal
and evaluate the signal at these points to obtain
\begin{equation}
\label{5.5}
S_0(d_j)=\epsilon f(\frac{(2j+1)\pi}{2\omega}). 
\end{equation}
Computing the spline curve $S_{mid}$ for the data points $(d_j,S_0(d_j))$, 
and subtracting this from the original signal we obtain after one iteration
of this algorithm that the new signal is given by
\begin{equation}
\label{5.6}
S_1^n(t)=S_0(t) - S_{mid}(t). 
\end{equation}
To compare the noise reduction efficiency of the two algorithms for this signal
on a finite time interval $[0,\frac{(2n+1)\pi}{\omega}]$ 
(i.e $k=0,\ldots,n$ and $j=0,\ldots,2n$)
we project the new signals on $\{\cos(\omega t),\sin(\omega t)\}$. 
(Both $\{\cos(\omega t),\sin(\omega t)\}$ have to be considered due to 
a possible phase shift in the new signal). To this end we have to compute
\begin{equation}
\label{5.7}
P_1^c=\displaystyle\int_{q_0}^{p_n}S_1^c(t)\cos(\omega t) dt,\,\,\
P_2^c=\displaystyle\int_{q_0}^{p_n}S_1^c(t)\sin(\omega t)dt.
\end{equation} 
and
\begin{equation}
\label{5.8}
Q_1^n=\displaystyle\int_{d_0}^{d_{2n}}S_1^n(t)\cos(\omega t) dt,\,\,\
Q_2^n=\displaystyle\int_{d_0}^{d_{2n}}S_1^n(t)\sin(\omega t)dt.
\end{equation} 

Using (\ref{5.1})-(\ref{5.6}) yields
\begin{equation}
\label{5.9}
P_1^c=\frac{\pi (k-1)}{2\omega}+
\displaystyle\int_{q_0}^{p_n}[\epsilon f(t)-\frac{S_{max}(t)+S_{min}(t)}{2}]\cos(\omega t)dt
\end{equation}
\begin{equation}
\label{5.10}
P_2^c= \displaystyle\int_{q_0}^{p_n}[\epsilon f(t)-
\frac{S_{max}(t)+S_{min}(t)}{2}]\sin(\omega t)dt
\end{equation} 
\begin{equation}
\label{5.11}
Q_1^n=\frac{\pi k}{2\omega}+
\displaystyle\int_{d_0}^{d_{2n}}[\epsilon f(t)-S_{mid}(t)]\cos(\omega t)dt
\end{equation}
\begin{equation}
\label{5.12}
Q_2^n= \displaystyle\int_{d_0}^{d_{2n}}[\epsilon f(t)-
S_{mid}(t)]\sin(\omega t)dt.
\end{equation} 
We conclude then that the efficiency of the algorithm to eliminate the noise
in the signal can be measured by the smallness of the absolute values of 
the integrals
$$
P_{mn}=\displaystyle\int_{q_0}^{p_n}[\epsilon f(t)-\frac{S_{max}(t)+S_{min}(t)}{2}]\cos(\omega t)dt,\,\,\,\,
Q_{mid}=\displaystyle\int_{d_0}^{d_{2n}}[\epsilon f(t)-S_{mid}(t)]\cos(\omega t)dt
$$
and the absolute values of $P_2^c$, $Q_2^n$.

To obtain a quantitative insight into this issue we considered the special 
case where
$$
f(t)= \cos(\nu t).
$$ 
with $k=0,\ldots,9$ and $j=0,\ldots,18$. A calculation of the $P_{mn}$ 
and the other integrals for $\nu \approx \omega$ yields:
$$
P_{mn}=\frac{26.703\epsilon}{\omega}+\frac{13.352\epsilon}{\omega^2}
(\nu-\omega) + O((\nu-\omega)^2)
$$
$$
Q_{mid}=\frac{28.274\epsilon}{\omega}+\frac{42.41\epsilon}{\omega^2}
(\nu-\omega) + O((\nu-\omega)^2)
$$
$$
P_2^c=\frac{-796.976\epsilon}{\omega^2}(\nu-\omega) +O((\nu-\omega)^2)
$$
$$
Q_2^n=\frac{-12.207\epsilon}{\omega^2}(\nu-\omega) +O((\nu-\omega)^2)
$$
These results show that when the frequency of the noise is close to the 
original frequency the classical algorithm leads to a large phase shift
in the signal and the noise is shifted with it.

In a more general setting of this analysis one may consider a Fourier 
expansion of $f(t)$ if this function is periodic.

For the convergence of the sifting iteration we now prove the following:

{\bf Theorem 3}: For the signal (\ref{5.1}) if we replace
the spline approximation between the midpoints by a linear interpolating 
function then the  sifting process will converge to $\cos(\omega t)$
if the derivatives of $f(t)$ in the $L^1$ norm are bounded. 

{\bf Proof}: The coordinates of the the midpoints between the max-min of 
(\ref{5.1}) are
$$
d_j=(\frac{(2j+1)\pi}{2\omega},\epsilon f(\frac{(2j+1)\pi}{2\omega})),
j=1,\ldots,k+1
$$
For a linear interpolation function eq. (\ref{5.6}) becomes
\begin{equation}
\label{5.28}
S_1^n(t)=\cos(\omega t)+\epsilon\left\{f(t)-\displaystyle\sum_{j=0}^k\frac{\omega\left[f(\frac{(2j+3)\pi}{2\omega})-f(\frac{(2j+1)\pi}{2\omega})\right](t-\frac{(2j+1)\pi}{2\omega})}{\pi}-f(\frac{(2j+1)\pi}{2\omega})\right\}
\end{equation}

The $L^1$ norm of the signal $S_0(t)$ is $O(1)$. To obtain an approximation 
for the $L^1$ norm of the perturbation $P_1=\epsilon f(t) -S_{mid}^1$
after one iteration we use trapezoidal integration. In this setting the 
integral of $f(t)$ cancels the integral of the linear interpolating function. 
This yields the following standard estimate for the residue of the perturbation
$$
\|P_1\|=\epsilon O((\frac{\pi}{\omega})^2)\|f^{\prime\prime}(t)\|
$$
(where primes denote differentiation with respect to $t$). We conclude that 
if the $L^1$ norm of derivatives of $f(t)$ are bounded then the sifting 
iterations will converge.

Using the same settings as in $theorem \,3$ (i.e replacing the
spline interpolating function by a linear interpolating function)
similar results apply to $Q_1^n$ and $Q_2^n$

{\bf Lemma\,7}: In the $L^1$ norm 
$$
\|Q_1^n-\frac{\pi k}{2\omega}\| =\epsilon O((\frac{\pi}{\omega})^2)\|f^{\prime\prime}(t)\|,\,\,\,\,
\|Q_1^n\| =\epsilon O((\frac{\pi}{\omega})^2)\|f^{\prime\prime}(t)\|
$$
The proof is the same as in $theorem\, 3$.

\subsection{Some Additional Analytical Insights} 

To obtain analytical insights about the performance of the 
EMD-midpoint algorithm we considered a signal of the form
\begin{equation}
\label{3.1a}
f(t)=\frac{1}{2}[\cos(\omega_1 t)+\cos(\omega_2t)],\,\,\,
\end{equation} 
where the ratio of the frequencies $\omega_1$, $\omega_2$ is a rational
number viz.
$$
\frac{\omega_2}{\omega_1}=\frac{m}{n}
$$ 
where $m,n$ are relative prime integers. In this case the signal
$f(t)$ is actually periodic with period $p=\frac{2n\pi}{\omega_1}$.
Due to this fact behavior of the classical versus the mid-point algorithm 
can be delineated without the need to discretize the signal.

On the interval $[0,p]$ the extrema of the signal which satisfy
$\frac{df}{dt}=0$ are given by
$$
\frac{\sin\omega_1 t}{\sin\omega_2 t}=-\frac{\omega_2}{\omega_1}=-\frac{m}{n}
$$
Computing these extrema points it is straightforward to construct the
spline approximations $S_{max}(t)$, $S_{min}(t)$ to the maximum and
minimum points and compute their average. Similarly we can
find the midpoints between the maxima and minima and evaluate the
corresponding spline approximation $S_{mid}(t)$ to the signal at 
these points. After one iteration of the sifting process the "sifted signal"
is given respectively by
\begin{equation}
\label{3.2}
h_{mn}(t)=f(t)-\frac{S_{max}(t)+S_{min}(t)}{2},
\end{equation} 
and
\begin{equation}
\label{3.3}
h_{mid}(t)=f(t)-S_{mid}(t).
\end{equation} 
The efficiency of the two algorithms can be deduced by projecting
these new signals on the Fourier components of the original signal. To
this end we compute
\begin{equation}
\label{3.4}
a_{mn}=\displaystyle\int_{0}^{p}h_{mn}(t)\cos(\omega_4 t) dt,\,\,\
b_{mn}=\displaystyle\int_{0}^{p}h_{mn}(t)\sin(\omega_4 t)dt.
\end{equation} 
\begin{equation}
\label{3.5}
c_{mn}=\displaystyle\int_{0}^{p}h_{mn}(t)\cos(\omega_5 t) dt,\,\,\
d_{mn}=\displaystyle\int_{0}^{p}h_{mn}(t)\sin(\omega_5 t)dt.
\end{equation} 
and
\begin{equation}
\label{3.6}
a_{mid}=\displaystyle\int_{0}^{p}h_{mid}(t)\cos(\omega_4 t) dt,\,\,\
b_{mid}=\displaystyle\int_{0}^{p}h_{mid}(t)\sin(\omega_5 t)dt.
\end{equation} 
\begin{equation}
\label{3.7}
c_{mid}=\displaystyle\int_{0}^{p}h_{mid}(t)\cos(\omega_4 t) dt,\,\,\
d_{mid}=\displaystyle\int_{0}^{p}h_{mid}(t)\sin(\omega_5 t)dt.
\end{equation} 
The amplitude of the Fourier components of the two frequencies in the 
classical EMD algorithm is 
\begin{equation}
\label{3.8}
A_{mn}=\sqrt{a_{mn}^2+b_{mn}^2},\,\,\,\, B_{mn}=\sqrt{c_{mn}^2+d_{mn}^2}.
\end{equation} 
Similarly for the mid-point algorithm we
\begin{equation}
\label{3.9}
A_{mid}=\sqrt{a_{mid}^2+b_{mid}^2},\,\,\,\, B_{mid}=\sqrt{c_{mid}^2+d_{mid}^2}.
\end{equation} 
The objective of the sifting process is to eliminate one of the
Fourier components in favor of the other. As a result the first IMF 
will contains, upon convergence, only one of the Fourier components in the 
original signal. Therefore the efficiency of the two algorithm can be 
inferred by comparing $A_{mn}$ versus $B_{mn}$ and $A_{mid}$ versus
$B_{mid}$. 

In the particular case where the signal is given by
\begin{equation}
\label{3.1}
f(t)=\frac{1}{2}[\cos(\omega_1 t)+\cos(\omega_2 t)],\,\,\,
\omega_1=\frac{3\pi}{64},\,\,\,\omega_2=\frac{\pi}{32}.
\end{equation}
$p=128$ (See Fig. $5$). Computing the integrals that appear in eqs.
(\ref{3.4})-(\ref{3.7}) we obtain
\begin{equation}
\label{3.11}
A_{mn}=31.63346911,\,\,\,\, B_{mn}=29.70292046,
\end{equation}
\begin{equation}
\label{3.12}
A_{mid}=34.19647843,\,\,\,\, B_{mid}=20.81145369.
\end{equation} 
These results show that after one iteration the classical EMD did not 
separate the two frequencies effectively ($A_{mn}$ and $B_{mn}$ are close 
to each other). On the other hand the mid-point algorithm performed well.

\setcounter{equation}{0}
\section{Examples and Comparisons}

Extensive numerical experiments were made to test and verify the efficiency of 
the modified algorithm. We present here the results of one of these tests 
in which the signal contains three close frequencies (where the classical
EMD algorithm fails). In our tests we considered also the effects of noise 
and phase shifts among the different frequencies but these will not presented 
here.
\begin{equation}
\label{2.1}
f(t)= \frac{1}{3}[\cos(\omega_1 t)+\cos(\omega_2 t)+\cos(\omega_3 t)]
\end{equation} 
where
$$
\omega_1=12\omega_0,\,\,\, \omega_2=10\omega_0,\,\,\, \omega_3=8\omega_0,
\,\,\, \omega_0=\frac{\pi}{256}.
$$
To apply the new EMD algorithm to this signal, discretized 
it over the time interval $[-2048,2048]$ by letting 
$t_{k+1}-t_{k}=1$, $k=1,\ldots,4097$.

The results of the signal decompositions into IMFs are presented in 
figures $1-4$. In all these figures the red lines represent the frequencies 
in the original signal (or its power spectrum) and the blue lines the 
corresponding intrinsic mode functions or their power spectrum which were 
obtained by the midpoint algorithm. 

Fig. $1$ is a plot of the data for the signal described by (\ref{2.1}).
Fig. $2$ represents the first IMF in the decomposition (versus 
the leading frequency in the data) while Figs. $3-4$ depict the spectral 
density distribution for the first two IMFs versus those related to 
the original frequencies in the data. It should be observed that although 
the amplitude of the spectral densities in these plots are somewhat different 
the maxima of the spectral density in each plot is very close to the 
original one.

\subsection{Cubic Lagrange Interpolation}

In both classical and the new versions of the EMD algorithm splines
are used for interpolation purposes. However the coefficients of each
spline polynomial depend non-locally on the data. As a result these 
coefficients might change if additional data is added. To compare
the two algorithm without this non-local dependence we replaced
the spline interpolation by cubic Lagrange interpolation (where 
the coefficients of the interpolating polynomial depend on the local
values of the function on the interval).

To carry out this comparison between the two EMD algorithms we considered
a signal composed of two frequencies and noise,
\begin{equation}
\label{5.13}
S_0(t)=\cos(\omega t) +\cos(\frac{3}{2}\omega t)+\epsilon f(t),\,\,\,\,
\epsilon \ll 1,
\end{equation}
on the time interval $[0,\frac{5\pi}{\omega}]$. 
The time interval was chosen so that the signal (without the perturbation)
has four maxima and minima on this interval. 

As in subsection $2.1$, we assume that the locations of the maxima
and minima do not change appreciably due to the perturbation. These 
locations are then given respectively by 
\begin{equation}
\label{5.14}
p_0=0,\,\,\, p_1=\frac{2}{\omega}\left(\pi-arctan\left(\frac{\sqrt{25-2\sqrt{10}}}{1+\sqrt{10}}\right)\right),\,\,\, p_2=\frac{4\pi}{\omega}-p_1,\,\,\,p3=\frac{4\pi}{\omega},
\end{equation}
\begin{equation}
\label{5.15}
q_0=\frac{2}{\omega}\left(arctan\left(\frac{\sqrt{25+2\sqrt{10}}}{-1+\sqrt{10}}\right)\right),\,\,\, q_1=\frac{4\pi}{\omega},\,\,\,q_2=\frac{4\pi}{\omega}-q_0,\,\,\,\, q_3=\frac{4\pi}{\omega}+q_0.
\end{equation}

Computing the cubic Lagrange interpolating polynomials $L_{max}$ and 
$L_{min}$ for the maxima and minima respectively we obtain for the 
"modified Lagrange classical-EMD" after one iteration
\begin{equation}
\label{5.16}
S_1^c(t)=S_0(t) - \frac{L_{max}(t)+L_{min}(t)}{2}.
\end{equation}
The number of midpoints on the interval $[0,\frac{5\pi}{\omega}]$ is seven.
For this reason we use two cubic Lagrange interpolating polynomials on
this interval. (The first is valid over the interval $[d_0,d_3]$ and 
the second is valid over $[d_3,d_6]$). Denoting this combined polynomial by
$L_{mid}$ and subtracting from the original signal we obtain 
after one iteration of this algorithm that the new signal is given by
\begin{equation}
\label{5.17}
S_1^n(t)=S_0(t) - L_{mid}(t). 
\end{equation}

To examine the performance of the two algorithms we project these new 
signals on $\cos(\frac{3}{2}\omega)$ and $\sin(\frac{3}{2}\omega)$.
\begin{equation}
\label{5.18}
P_1^c=\displaystyle\int_{q_0}^{p_3}S_1^c(t)\cos(\frac{3}{2}\omega t) dt,\,\,\
P_2^c=\displaystyle\int_{q_0}^{p_3}S_1^c(t)\sin(\frac{3}{2}\omega t)dt.
\end{equation}
and
\begin{equation}
\label{5.19}
Q_1^n=\displaystyle\int_{d_0}^{d_{6}}S_1^n(t)\cos(\frac{3}{2}\omega t) dt,\,\,\
Q_2^n=\displaystyle\int_{d_0}^{d_{6}}S_1^n(t)\sin(\frac{3}{2}\omega t)dt.
\end{equation}
Furthermore if we assume that $f=\cos(\nu t)$ and
$\nu \approx \frac{3}{2}\omega$ we obtain to order $\epsilon$
\begin{equation}
\label{5.22}
P_1^c=\frac{3.8568}{\omega}+\frac{0.0175\epsilon}{\omega}+
O(\nu-\frac{3}{2}\omega),\,\,\,\,
P_2^c=-\frac{1.0637}{\omega}+\frac{0.0399\epsilon}{\omega}+
O(\nu-\frac{3}{2}\omega)
\end{equation}
\begin{equation}
\label{5.23}
Q_1^n=\frac{6.3795}{\omega}-\frac{0.1257\epsilon}{\omega}+
O(\nu-\frac{3}{2}\omega),\,\,\,\,
Q_2^n=-\frac{0.2184}{\omega}+\frac{0.3113\epsilon}{\omega}+
O(\nu-\frac{3}{2}\omega).
\end{equation}
These results demonstrate the superiority on the midpoint algorithm 
in this setting.(The total projection of the new signal on 
$\cos\frac{3}{2}\omega t$ is larger and the phase shift is smaller).

\setcounter{equation}{0}
\section{Convergence Rates}

To compare the convergence rates of the classical versus the midpoint 
algorithm we considered three cases all of which were composed of two 
frequencies. In the first case the two frequencies were well separated.
In the second case the two frequencies were close while in the third case
they were almost "overlapping". In all cases the signal was given by
$$
f(t)=\frac{1}{2}(\cos\omega_1 t +\cos\omega_2 t)
$$
This signal was discretized on the time interval $[-2048,2048]$ 
with $\Delta t=1$.
 
For the first case the two frequencies were
$$
\omega_1=12\omega,\,\,\,\, \omega_2=8\omega,\,\,\,\, \omega=\frac{\pi}{256}.
$$
As can be expected both the classical and midpoint algorithm were able to
discern the individual frequencies through the sifting algorithm. However 
it took the classical algorithm $59$ iterations to converge to the first IMF. 
On the other hand the midpoint algorithm converged in only $7$ iterations 
(using the same convergence criteria). We wish to point out also that the
midpoint algorithm has a lower computational cost than the classical 
algorithm. It requires in each iteration the computation of only one spline 
interpolating polynomial. On the other hand the classical algorithm requires 
two such polynomials, one for the maximum points and one for the minimum 
points.

For the second test the frequencies were
$$
\omega_1=\frac{\pi}{24}+\frac{\pi}{288},\,\,\,\,
\omega_2=\frac{\pi}{24}-\frac{\pi}{288}
$$
that is the difference between the two frequencies is $\frac{\pi}{144}$.

In this case the midpoint algorithm was able to separate the two frequencies.
Fig $6$ and Fig $7$ compare the power spectrum of the original frequencies 
versus those of $IMF_1$ and $IMF_2$ which were obtained through this algorithm. 
Convergence to $IMF_1$ was  obtained in 18 iterations and $IMF_2$ was obtained 
by $7$ additional iterations.

The classical EMD algorithm did converge to $IMF_1$ in $45$ iterations but 
the power spectrum of this $IMF$ deviated significantly from the first 
frequency in the signal. $IMF_2$ failed (completely) to 
detect  correctly the second frequency.

In third case the frequencies were
$$
\omega_1=\frac{\pi}{24}+\frac{\pi}{1000},\,\,\,\,
\omega_1=\frac{\pi}{24}-\frac{\pi}{1000}.
$$
In this case the classical algorithm was unable to separate the two 
frequencies i.e $IMF_1$ contained both frequencies. 
The midpoint algorithm did somewhat better but the resolution was not 
complete. Moreover the sifting process in both cases led 
to the creation of "ghost frequencies" which were not present in the 
original signal. 

At this juncture one might wonder if a "hybrid algorithm" whereby the sifting
function is the average (or some similar combination) of those obtained by 
the classical and midpoint algorithms might outperform the separate algorithms
(in spite of the obvious additional computational cost). However our 
experimentations with such algorithm did not yield the desired results
(i.e. the convergence rate and resolution did not improve).

\section{Conclusion}

In this paper we presented a variant of the EMD algorithm which utilizes 
the midpoints between the max-min points of the signal in the sifting
iterative process. We demonstrated through several case studies and 
theoretical approximations that this algorithm can resolve signals with 
moderately close frequencies where the classical EMD algorithm fails. 
We showed also that it has a better convergence rate. From a formal point 
of view this superior performance of the midpoint algorithm can be traced 
to the fact that the deviation of the signal average from zero is sampled 
at "half" the scale of the classical EMD algorithm.

\section*{References}

\begin{itemize}

\item[1] N. E. Huang - USA Patent $\# 6,311,130 B 1$ , Date Oct 30,2001

\item[2] N. E. Huang et all, ``The empirical mode decomposition and the Hilbert
spectrum for nonlinear and non-stationary time series analysis",
Proceedings of the Royal Society Vol. 454 pp. 903-995 (1998)

\item[3] Gabriel Rilling and Patrick Flandrin,  ``One or Two Frequencies? The
Empirical Mode Decomposition Answers", IEEE Trans. Signal Analysis 
Vol. 56 pp.85-95 (2008).

\item[4] Zhaohua Wu and Norden E. Huang, ``On the Filtering Properties of 
the Empirical Mode Decomposition, Advances in Adaptive Data Analysis", 
Volume: 2, Issue: 4 pp. 397-414. (2010)

\item[5] Albert Ayenu-Prah and  Nii Attoh-Okine, ``A Criterion for Selecting 
Relevant Intrinsic Mode Functions in Empirical Mode Decomposition",
Advances in Adaptive Data Analysis, Vol. 2, Issue: 1(2010) pp. 1-24.

\item[6] G. Rilling, P. Flandrin and P. Goncalves, ``Empirical 
Mode Decomposition As a Filter Bank, IEEE Signal Processing Letters,
vol. 11, no. 2, pp. 112-114, 2004

\end{itemize}

\newpage
\begin{figure}[ht!]
\includegraphics[scale=1,height=160mm,angle=0,width=180mm]{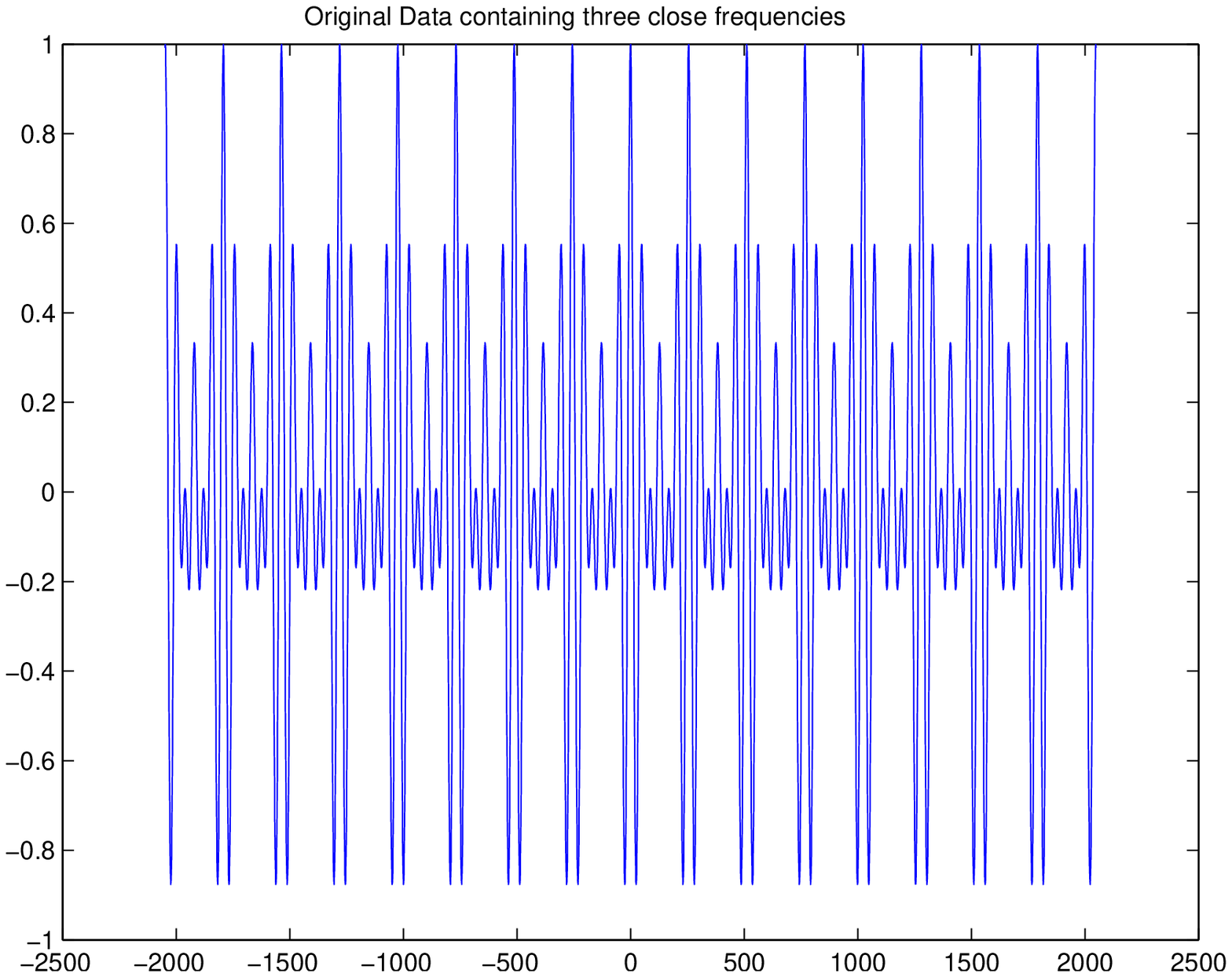}
\caption{}
\end{figure}
\newpage
\begin{figure}[ht!]
\includegraphics[scale=1,height=160mm,angle=0,width=180mm]{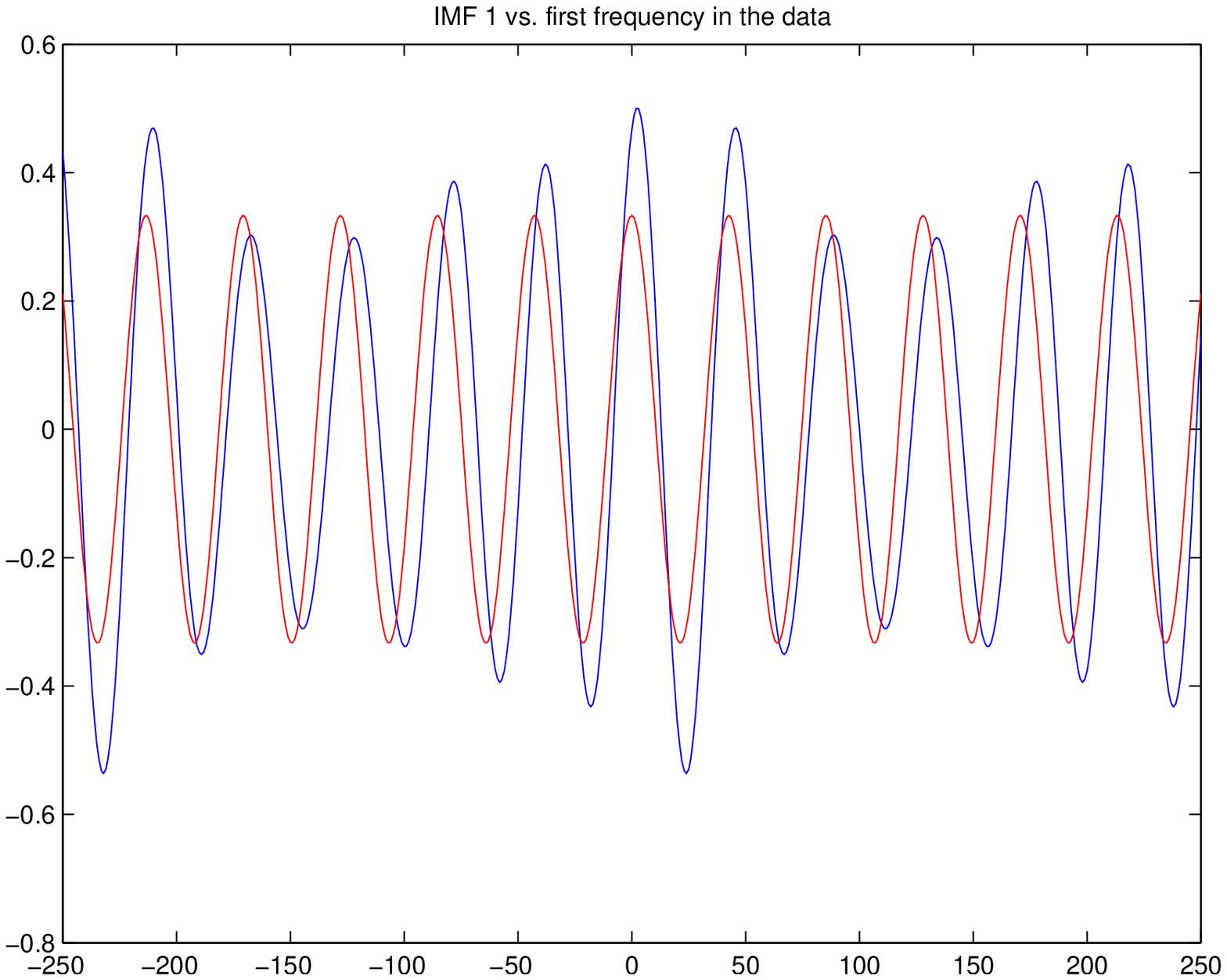}
\caption{}
\end{figure}
\newpage
\begin{figure}[ht!]
\includegraphics[scale=1,height=160mm,angle=0,width=180mm]{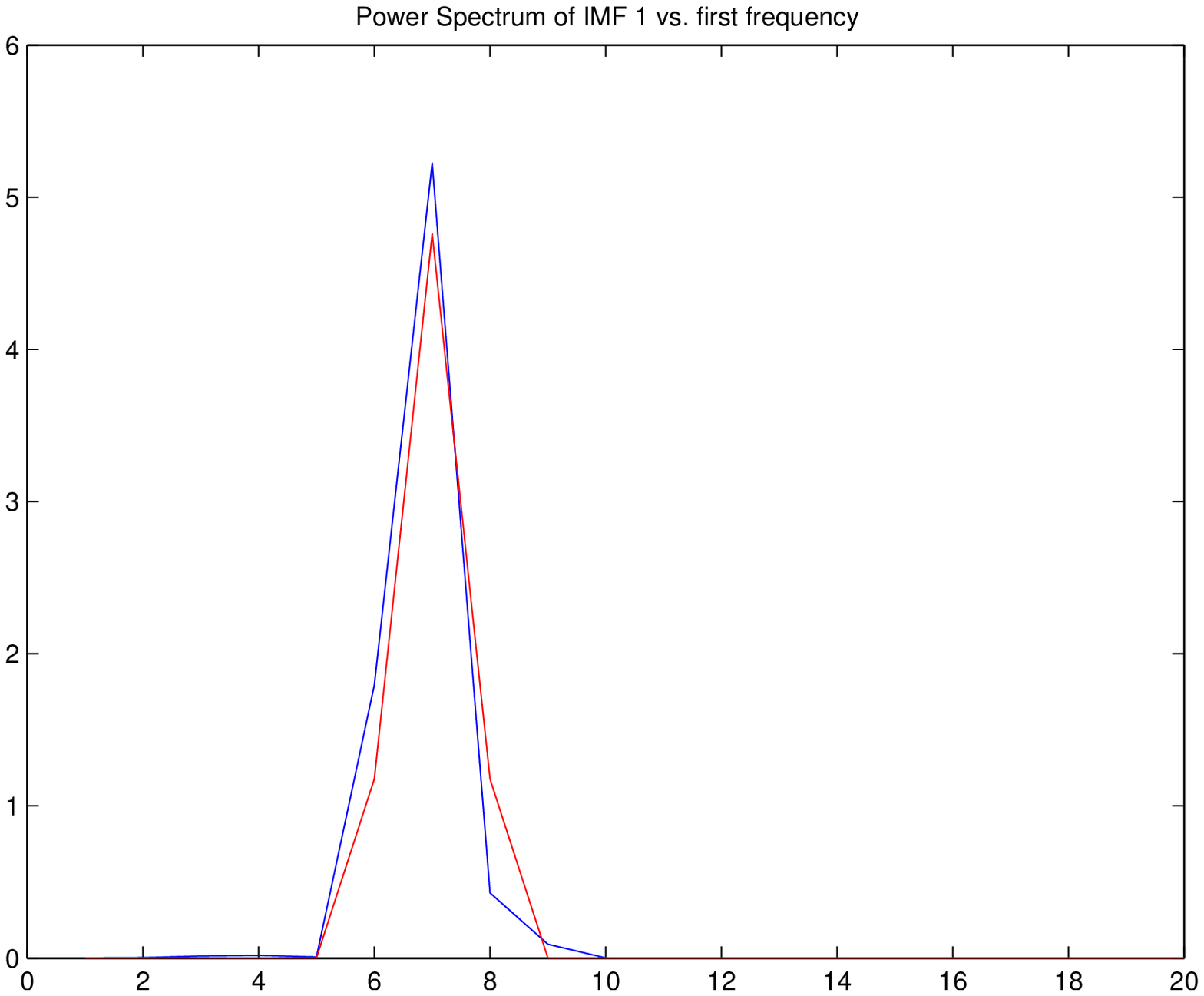}
\caption{}
\end{figure}
\newpage
\begin{figure}[ht!]
\includegraphics[scale=1,height=160mm,angle=0,width=180mm]{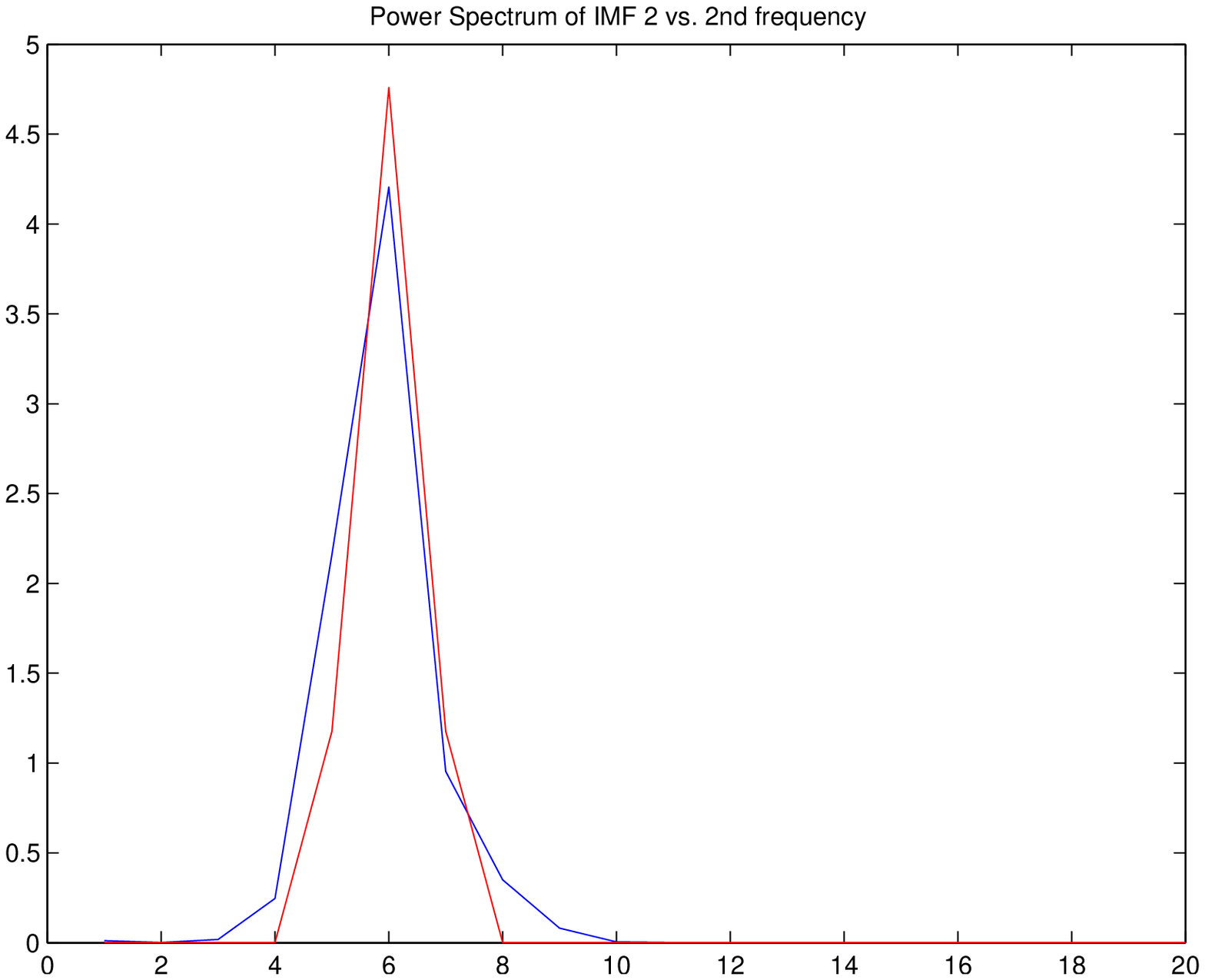}
\caption{}
\end{figure}
\newpage
\begin{figure}[ht!]
\includegraphics[scale=1,height=160mm,angle=0,width=180mm]{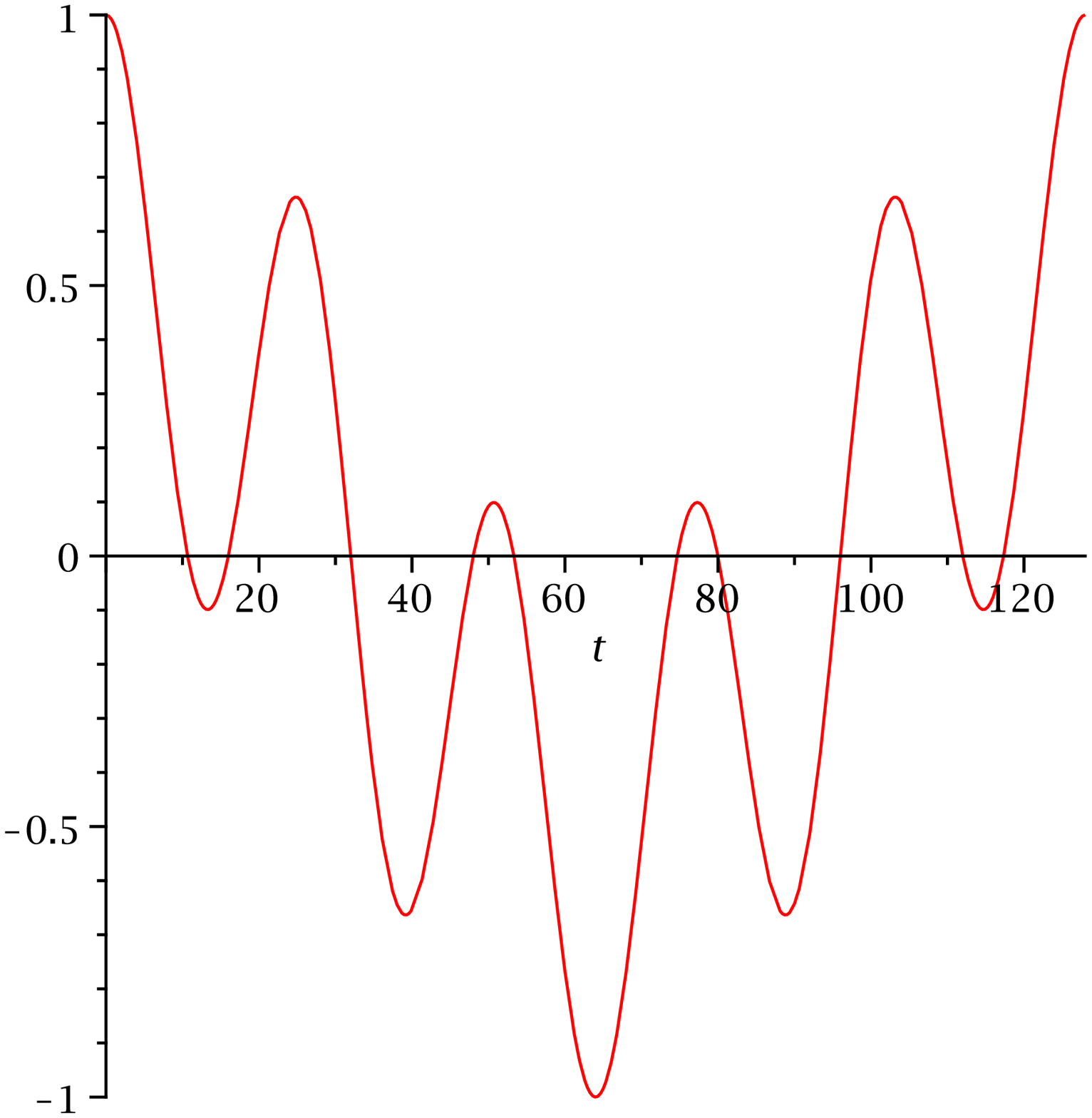}
\caption{}
\end{figure}
\newpage
\begin{figure}[ht!]
\includegraphics[scale=1,height=160mm,angle=0,width=180mm]{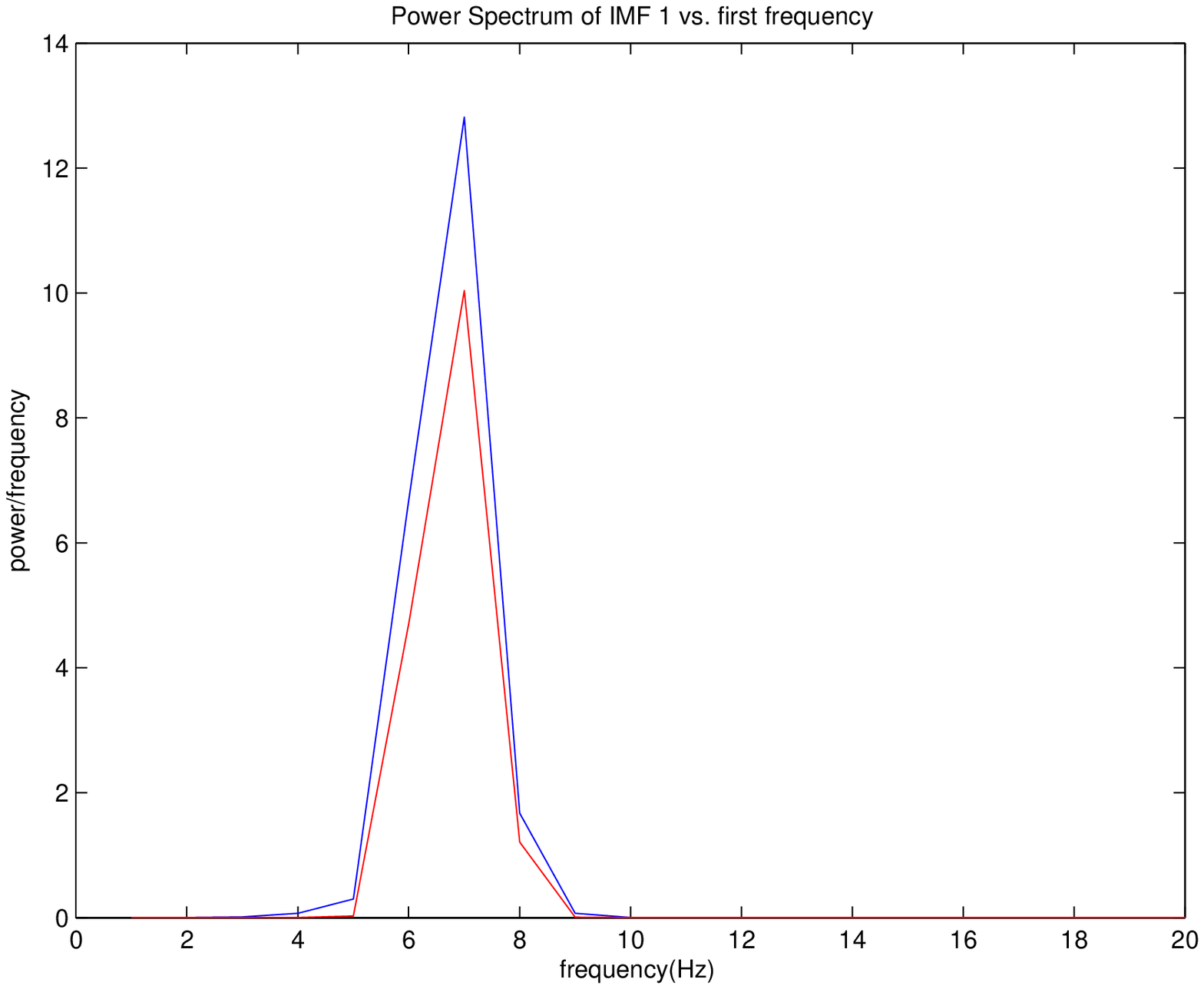}
\caption{}
\end{figure}
\newpage
\begin{figure}[ht!]
\includegraphics[scale=1,height=160mm,angle=0,width=180mm]{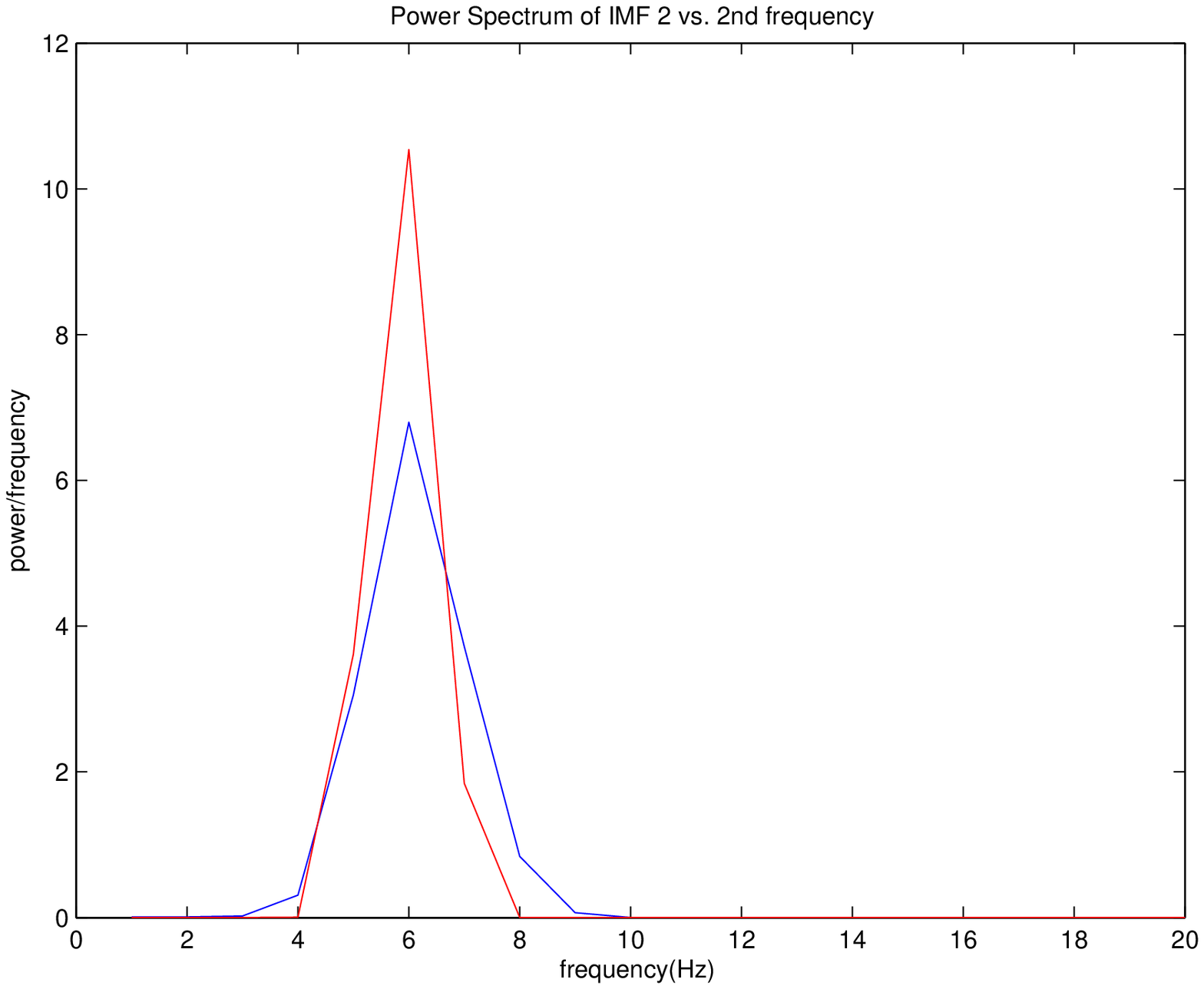}
\caption{}
\end{figure}
\end{document}